# NONPARAMETRIC ESTIMATION BY CONVEX PROGRAMMING

By Anatoli B. Juditsky and Arkadi S. Nemirovski[1]

*Université Grenoble I and Georgia Institute of Technology*

The problem we concentrate on is as follows: given (1) a convex compact set $X$ in $\mathbb{R}^n$, an affine mapping $x \mapsto A(x)$, a parametric family $\{p_\mu(\cdot)\}$ of probability densities and (2) $N$ i.i.d. observations of the random variable $\omega$, distributed with the density $p_{A(x)}(\cdot)$ for some (unknown) $x \in X$, estimate the value $g^T x$ of a given linear form at $x$.

For several families $\{p_\mu(\cdot)\}$ with no additional assumptions on $X$ and $A$, we develop computationally efficient estimation routines which are minimax optimal, within an absolute constant factor. We then apply these routines to recovering $x$ itself in the Euclidean norm.

**1. Introduction.** The problem we are interested in is essentially as follows: suppose that we are given a convex compact set $X$ in $\mathbb{R}^n$, an affine mapping $x \mapsto A(x)$ and a parametric family $\{p_\mu(\cdot)\}$ of probability densities. Suppose that $N$ i.i.d. observations of the random variable $\omega$, distributed with the density $p_{A(x)}(\cdot)$ for some (unknown) $x \in X$, are available. Our objective is to estimate the value $g^T x$ of a given linear form at $x$.

In nonparametric statistics, there exists an immense literature on various versions of this problem (see, e.g., [10, 11, 12, 13, 15, 17, 18, 21, 22, 23, 24, 25, 26, 27, 28] and the references therein). To the best of our knowledge, the majority of papers on the subject focus on specific domains $X$ (e.g., distributions with densities from Sobolev balls), and investigate lower and upper bounds on the worst-case, with regard to $x \in X$, accuracy to which the problem of interest can be solved. These bounds depend on the number of observations $N$, and the question of primary interest is the behavior of those bounds as $N \to \infty$. When the lower and the upper bounds coincide within a constant factor [or, ideally, within factor $(1 + o(1))$ as $N \to \infty$], the estimation problem is considered essentially solved, and the estimation methods underlying the upper bounds are treated as optimal.

Received March 2008; revised July 2008.
[1]Supported in part by the NSF Grant 0619977.

The approach we adopt in this paper is of a different spirit; we make no "structural assumptions" on $X$, aside from assumptions of convexity and compactness which are crucial for us, and we make no assumptions on the linear functional $p$. Clearly, with no structural assumptions on $X$ and $p$, explicit bounds on the risks of our estimates, as well as bounds on the minimax optimal risk, are impossible. However, it is possible to show that *when estimating linear forms, the worst-case risk of the estimator we propose is within an absolute constant factor of the "ideal"* (i.e., the minimax optimal) *risk*. It should be added that while the optimal, within an absolute constant factor, worst-case risk of our estimates is not available in a closed analytical form, it is "available algorithmically"—it can be efficiently computed, provided that $X$ is computationally tractable.[1]

Note that the estimation problem, presented above, can be seen as a generalization of the problem of estimation of linear functionals of the central parameter of a normal distribution (see [4, 8, 9, 16]). Namely, suppose that the observation $\omega \in \mathbb{R}^m$,

$$\omega = Ax + \sigma\xi$$

of the unknown signal $x$ is available. Here $A$ is a given $m \times n$ matrix and $\xi \sim \mathcal{N}(0, I_m)$, $\sigma > 0$ is known. For this important case the problem has been essentially solved in [5], where it was proved that for several commonly used loss functions, the minimax optimal *affine in $\omega$* estimate is minimax optimal, within an absolute constant factor, among *all possible* estimates.

Another special case of our setting is the problem of estimating a linear functional $g(p)$ of an unknown distribution $p$, given $N$ i.i.d. observations $\omega_1, \ldots, \omega_N$, which obey $p$. We suppose that it is known a priori that $p \in X$, where $X$ is a given convex compact set of distributions (here the parameter $x$ is the density $p$ itself). Some important results for this problem have been obtained in [6] and [7]. For instance, in [7] the authors established minimax bounds for the risk of estimation of $g(p)$ and developed an estimation method based on the binary search algorithm. The estimation procedure uses at each search iteration tests of convex hypotheses, studied in [2, 3]. That estimator of $g(p)$ is shown to be minimax optimal (within an absolute constant factor) if some basic structural assumptions about $X$ hold.

In this paper, we concentrate on the properties of *affine estimators*. Here, we refer to an estimator $\hat{g}$ as affine when it is of the form $\hat{g}(\omega_1, \ldots, \omega_N) = \sum_{i=1}^{N} \phi(\omega_i)$, for some given functions $\phi$, that is, if $\hat{g}$ is an affine function of the empirical distribution. When $\phi$ itself is an affine function, the estimator is also affine in the observations, as it is in the setting of [5]. Our motivation is to extend the results obtained in [5] to the non-Gaussian situation. In particular, we propose a technique of derivation of affine estimators which are minimax optimal (up to a moderate absolute constant) for a class of "good parametric families of distributions," which is defined in Section 2.1. As normal family and discrete distributions belong to the class of good parametric families, the minimax optimal estimators for these cases are obtained by direct application of the general construction. In this sense, our results generalize those of [7] and [5] on the estimation of linear functionals. On the other hand, it is clear that different techniques, presented in the current paper, inherit from those developed in [3] and [7]. To make a computationally efficient solution of the estimation problem possible, unlike the authors of those papers, we concentrate only on the finite-dimensional situation. As a result, the proposed estimation procedures allow efficient numeric implementation. This also allows us to avoid much of the intricate mathematical details. However, we allow the dimension to be arbitrarily large, thus addressing, essentially, a nonparametric estimation problem.

---

[1] For details on computational tractability and complexity issues in convex optimization, see, for example, [1], Chapter 4. A reader not familiar with this area will not lose much when interpreting a computationally tractable convex set as a set given by a finite system of inequalities $p_i(x) \leq 0$, $i = 1, \ldots, m$, where $p_i(x)$ are convex polynomials.

The rest of this paper is organized as follows. In Section 2, we define the main components of our study—we state the estimation problem and define the corresponding risk measures. Then, in Section 3, we provide the general solution to the estimation problem, which is then applied, in Section 4, to the problems of estimating linear functionals in the normal model and the tomography model. Finally, in Section 5, we present adaptive versions of affine estimators.

Note that when passing from recovering linear forms of the unknown signal to recovering the signal itself, we do impose structural assumptions on $X$, but still make no structural assumptions on the affine mapping $A(x)$. Our "optimality results" become weaker—instead of "optimality within an absolute constant factor" we end up with statements like "the worst-case risk of such-and-such estimate is in between the minimax optimal risk and the latter risk to the power $\chi$," with $\chi$ depending on the geometry of $X$ (and close to 1 when this geometry is "good enough").

## 2. Problem statement.

2.1. *Good parametric families of distributions.* Let $(\Omega, P)$ be a Polish space with Borel $\sigma$-finite measure, and $\mathcal{M} \subset \mathbb{R}^m$. Assume that every $\mu \in \mathcal{M}$ is associated with a probability density $p_\mu(\omega)$—a Borel nonnegative function on $\Omega$ such that $\int_\Omega p_\mu(\omega) P(d\omega) = 1$; we refer to the mapping $\mu \to p_\mu(\cdot)$ as to a *parametric density family* $\mathcal{D}$. Let also $\mathcal{F}$ be a finite-dimensional linear space of Borel functions on $\Omega$ which contains constants. We call a pair $(\mathcal{D}, \mathcal{F})$ *good* if it possesses the following properties:

1. $\mathcal{M}$ is an open convex set in $\mathbb{R}^m$;
2. whenever $\mu \in \mathcal{M}$, we have $p_\mu(\omega) > 0$ everywhere on $\Omega$;
3. whenever $\mu, \nu \in \mathcal{M}$, we have $\phi(\omega) = \ln(p_\mu(\omega)/p_\nu(\omega)) \in \mathcal{F}$;
4. whenever $\phi(\omega) \in \mathcal{F}$, the function

$$F_\phi(\mu) = \ln\left(\int_\Omega \exp\{\phi(\omega)\} p_\mu(\omega) P(d\omega)\right)$$

is well defined and concave in $\mu \in \mathcal{M}$.

The reader familiar with exponential families will immediately recognize that the above definition implies that $\mathcal{D}$ is such a family. Let us denote $p_\mu(\omega) = \exp\{\theta(\mu)^T \omega - C(\theta(\mu))\}, \mu \in \mathcal{M}$, its density with regard to $P$ where $\theta$ is the natural parameter and $C(\cdot)$ as the cumulant function. Then, $\mathcal{D}$ is good if:

1. $\mathcal{M}$ is an open convex set in $\mathcal{D}_P = \{\mu \in \mathbb{R}^m | \int e^{\theta(\mu)^T \omega} P(d\omega) < \infty\}$;
2. for any $\phi$ such that the cumulant function $C(\theta(\mu) + \phi)$ is well definded, the function $[C(\theta(\mu) + \phi) - C(\theta(\mu))]$ is concave in $\mu \in \mathcal{M}$.

Let us list several examples.

EXAMPLE 1 (*Discrete distributions*). Let $\Omega = \{1, 2, \ldots, M\}$ be a finite set, $P$ be a counting measure on $\Omega$, $\mathcal{M} = \{\mu \in \mathbb{R}^M : \mu > 0, \sum_i \mu_i = 1\}$ and $p_\mu(i) = \mu_i$, $i = 1, \ldots, M$. Let also $\mathcal{F}$ be the set of all functions on $\Omega$. The associated pair $(\mathcal{D}, \mathcal{F})$ clearly is good.

EXAMPLE 2 (*Poisson distributions*). Let $\Omega = \{0, 1, \ldots\}$, $P$ be the counting measure on $\Omega$, $\mathcal{M} = \{\mu \in \mathbb{R} : \mu > 0\}$ and $p_\mu(i) = \frac{\mu^i \exp\{-\mu\}}{i!}$, $i \in \Omega$, so that $p_\mu$ is the Poisson distribution with the parameter $\mu$. Let also $\mathcal{F}$ be the set of affine functions $\phi(i) = \alpha i + \beta$ on $\Omega$. We claim that the associated pair $(\mathcal{D}, \mathcal{F})$ is good. Indeed, $\ln(p_\mu(i)/p_\nu(i)) = i[\ln \mu - \ln \nu] + \mu - \nu$ is an affine function of $i$, and

$$\ln\left(\sum_i \exp\{\alpha i + \beta\} \frac{\mu^i \exp\{-\mu\}}{i!}\right) = \ln(\exp\{\beta - \mu\} \exp\{\mu \exp\{\alpha\}\})$$
$$= \beta - \mu + \mu \exp\{\alpha\}$$

is a concave function of $\mu > 0$.

EXAMPLE 3 (*Gaussian distributions with fixed covariance*). Let $\Omega = \mathbb{R}^k$, $P$ be the Lebesque measure on $\Omega$, $\Sigma$ be a positive definite $k \times k$ matrix, $\mathcal{M} = \mathbb{R}^k$ and

$$p_\mu(\omega) = (2\pi)^{-k/2}(\mathrm{Det}\,\Sigma)^{-1/2} \exp\{-(\omega - \mu)^T \Sigma^{-1}(\omega - \mu)\}$$

be the density of the Gaussian distribution with mean $\mu$ and covariance matrix $\Sigma$. Let, further, $\mathcal{F}$ be comprised of affine functions on $\Omega$. We claim that the associated pair $(\mathcal{D}, \mathcal{F})$ is good. Indeed, the function $\ln(p_\mu(\omega)/p_\nu(\omega))$ indeed is affine on $\Omega$, and

$$\ln\left(\int \exp\{\phi^T \omega + c\} p_\mu(\omega)\, d\omega\right) = c + \phi^T \mu + \frac{\phi^T \Sigma \phi}{2}$$

is a concave function of $\mu$.

EXAMPLE 4 (*Direct product of good pairs*). Let $p^\ell_{\mu_\ell}(\omega_\ell)$ be a probability density, parameterized by $\mu_\ell \in \mathcal{M}_\ell \subset \mathbb{R}^{m_\ell}$, on a Polish space $\Omega_\ell$ with Borel $\sigma$-finite measure $P_\ell$, and $\mathcal{F}_\ell$ be a finite-dimensional linear space of Borel functions on $\Omega_\ell$ such that the associated pairs $(\mathcal{D}_\ell, \mathcal{F}_\ell)$ are good. Let us define the *direct product* $(\mathcal{D}, \mathcal{F}) = \bigotimes_{\ell=1}^L (\mathcal{D}_\ell, \mathcal{F}_\ell)$ of these pairs as follows:

- The associated space with measure is $(\Omega = \Omega_1 \times \cdots \times \Omega_L, P = P_1 \times \cdots \times P_\ell)$.
- The set of parameters is $\mathcal{M} = \mathcal{M}_1 \times \cdots \times \mathcal{M}_L$, and the density associated with a parameter $\mu = (\mu_1, \ldots, \mu_L)$ from this set is $p_\mu(\omega_1, \ldots, \omega_L) = \prod_{\ell=1}^L p^\ell_{\mu_\ell}(\omega_\ell)$.
- $\mathcal{F}$ is comprised of all functions $\phi(\omega_1, \ldots, \omega_L) = \sum_{\ell=1}^L \phi_\ell(\omega_\ell)$ with $\phi_\ell(\cdot) \in \mathcal{F}_\ell$, $\ell = 1, \ldots, m$.

We claim that *the direct product of good pairs is good*. Indeed, $\mathcal{M}$ is an open convex set; when $\mu = (\mu_1, \ldots, \mu_L)$ and $\nu = (\nu_1, \ldots, \nu_L)$ are in $\mathcal{M}$, we have

$$\ln(p_\mu(\omega_1, \ldots, \omega_L)/p_\nu(\omega_1, \ldots, \omega_L)) = \sum_{\ell=1}^L \ln(p^\ell_{\mu_\ell}(\omega_\ell)/p^\ell_{\nu_\ell}(\omega_\ell)) \in \mathcal{F}$$

and when $\phi(\omega_1, \ldots, \omega_L) = \sum_\ell \phi_\ell(\omega_\ell) \in \mathcal{F}$, we have

$$\ln\left(\int_\Omega \exp\{\phi(\omega)\} p_\mu(\omega) P(d\omega)\right) = \ln\left(\prod_\ell \int_{\Omega_\ell} \exp\{\phi_\ell(\omega_\ell)\} p^\ell_{\mu_\ell}(\omega_\ell) P(d\omega_\ell)\right)$$

$$= \sum_\ell \ln\left(\int_{\Omega_\ell} \exp\{\phi_\ell(\omega_\ell)\} p^\ell_{\mu_\ell}(\omega_\ell) P(d\omega_\ell)\right),$$

which is a sum of concave functions of $\mu$ and thus is concave in $\mu$.

2.2. *The problem.* The problem we are interested in is as follows:

PROBLEM I. We are given the following:

- a convex compact set $X \subset \mathbb{R}^n$,
- a good pair $(\mathcal{D}, \mathcal{F})$ comprised of

– a parametric family $\{p_\mu(\omega) : \mu \in \mathcal{M} \subset \mathbb{R}^m\}$ of probability densities on a Borel space $\Omega$ with $\sigma$-finite Borel measure $P$ and
– a finite-dimensional linear space $\mathcal{F}$ of Borel functions on $\Omega$,
- an affine mapping $x \mapsto A(x) : X \mapsto \mathcal{M}$,
- a linear form $g^T z$ on $\mathbb{R}^n \supset X$.

Aside of this a priori information, we are given a realization $\omega$ of a random variable taking values in $\Omega$ and distributed with the density $p_{A(x)}(\cdot)$ for some *unknown in advance* $x \in X$. Our goal is to infer from this observation an estimate $\hat{g}(\omega)$ of the value $g^T x$ of the given linear form at $x$.

From now on we refer to an estimate as *affine*, if it is of the form $\hat{g}(\omega) = \phi(\omega)$, with certain $\phi \in \mathcal{F}$.

We quantify the risk of a candidate estimate $\hat{g}(\cdot)$ by its worst-case, over $x \in X$, confidence interval, given the confidence level. Specifically, given a *confidence level* $\varepsilon \in (0,1)$, we define the associated $\varepsilon$-*risk* of an estimate $\hat{g}$ as

$$\mathrm{Risk}(\hat{g}; \varepsilon) = \inf \left\{ \delta : \sup_{x \in X} \mathrm{Prob}_{\omega \sim p_{A(x)}(\cdot)} \{\omega : |\hat{g}(\omega) - g^T x| > \delta\} < \varepsilon \right\}.$$

The corresponding minimax optimal $\varepsilon$-risk is defined as

$$\mathrm{Risk}_*(\varepsilon) = \inf_{\hat{g}(\cdot)} \mathrm{Risk}(\hat{g}; \varepsilon),$$

where inf is taken over the space of all Borel functions $\hat{g}$ on $\Omega$. We are interested also in the minimax optimal $\varepsilon$-risk of affine estimates

$$\mathrm{RiskA}(\varepsilon) = \inf_{\phi(\cdot) \in \mathcal{F}} \mathrm{Risk}(\phi; \varepsilon).$$

## 3. Minimax optimal affine estimators.

### 3.1. *Main result.*  Our main result follows.

THEOREM 3.1. *Let the pair $(\mathcal{D}, \mathcal{F})$ underlying Problem I be good. Then, the minimax optimal risk achievable with affine estimates is, for small $\varepsilon$, within an absolute constant factor of the "true" minimax optimal risk, specifically,*

$$0 \leq \varepsilon < 1/4 \quad \Rightarrow \quad \mathrm{RiskA}(\varepsilon) \leq \theta(\varepsilon) \mathrm{Risk}_*(\varepsilon), \qquad \theta(\varepsilon) = \frac{2 \ln(2/\varepsilon)}{\ln(1/(4\varepsilon))}.$$

PROOF.  For $r \geq 0$, let us set

$$\Phi_r(x, y; \phi, \alpha) = g^T x - g^T y + \alpha \ln \left( \int_\Omega \exp\{\alpha^{-1} \phi(\omega)\} p_{A(y)}(\omega) P(d\omega) \right)$$

$$+ \alpha \ln\left(\int_\Omega \exp\{-\alpha^{-1}\phi(\omega)\} p_{A(x)}(\omega) P(d\omega)\right)$$
$$+ 2\alpha r: \qquad Z \times \mathcal{F}_+ \to \mathbb{R},$$
$$Z = X \times X,$$
$$\mathcal{F}_+ = \mathcal{F} \times \{\alpha > 0\}.$$

We claim that this function is a continuous real-valued function on $Z \times \mathcal{F}_+$, which is convex in $(\phi, \alpha) \in \mathcal{F}_+$ and concave in $(x,y) \in Z$.

Indeed, the function
$$\Psi(\mu,\nu;\phi) = \ln\left(\int_\Omega \exp\{\phi(\omega)\} p_\mu(\omega) P(d\omega)\right)$$
$$+ \ln\left(\int_\Omega \exp\{-\phi(\omega)\} p_\nu(\omega) P(d\omega)\right): \qquad (\mathcal{M} \times \mathcal{M}) \times \mathcal{F} \to \mathbb{R}$$

is well defined, concave in $(\mu,\nu) \in \mathcal{M} \times \mathcal{M}$ [since $(\mathcal{D}, \mathcal{F})$ is good] and convex in $\phi \in \mathcal{F}$ (evident). Since $\mathcal{M}$ is open and $\mathcal{F}$ is a finite-dimensional linear space, $\Psi$ is continuous on its domain. It remains to note that $\Phi_\varepsilon$ is the sum of a linear function of $x, y, \alpha$ and the function $\alpha \Psi(A(x), A(y); \alpha^{-1}\phi)$ which clearly is concave in $(x,y)$ [since $\Psi(\mu,\nu;\phi)$ is concave in $(\mu,\nu)$ and $A(\cdot)$ is affine] and convex in $(\phi, \alpha) \in \mathcal{F}_+$ [since $\Psi(\mu,\nu;\phi)$ is continuous in $\phi \in \mathcal{F}$, and the transformation $f(u) \mapsto g(u,\alpha) = \alpha f(u/\alpha)$ converts a convex function of $u$ into a convex in $(\alpha > 0, u)$ function of $(u, \alpha)$].

Since $Z$ is a convex finite-dimensional compact set, $\mathcal{F}_+$ is a convex finite-dimensional set and $\Phi_\varepsilon$ is continuous and convex–concave on $Z \times \mathcal{F}_+$, we can invoke the Sion–Kakutani theorem (see, e.g., [14]) to infer that

$$(3.1) \quad \sup_{x,y \in X} \inf_{\phi \in \mathcal{F}, \alpha > 0} \Phi_r(x,y;\phi,\alpha) = \inf_{\phi \in \mathcal{F}, \alpha > 0} \max_{x,y \in X} \Phi_r(x,y;\phi,\alpha) := 2\Phi_*(r).$$

Note that $\Phi_*(r) \geq 0$ is a concave and nonnegative function of $r \geq 0$. Indeed, the functional $f_x[h] = \ln \int_\Omega \exp\{h(\omega)\} p_{A(x)}(\omega) P(d\omega)$ is well defined and convex on $\mathcal{F}$, whence

$$\Phi_r(x,x;\phi,\alpha) = 2\alpha r + \alpha(f_x[-\alpha^{-1}\phi] + f_x[\alpha^{-1}\phi]) \geq 2\alpha r \geq 0,$$

whence $\Phi_*(r) \geq \frac{1}{2} \sup_{x \in X} \inf_{\phi \in \mathcal{F}, \alpha > 0} \Phi_r(x,x;\phi,\alpha) \geq 0$. The concavity of $\Phi_*(r)$ on the nonnegative ray follows immediately from the representation, yielded by (3.1),

$$\Phi_*(r) = \frac{1}{2} \inf_{\phi \in \mathcal{F}, \alpha}\left[2\alpha r + \sup_{x,y \in X} \Phi_0(x,y;\phi,\alpha)\right]$$

of $\Phi_*(r)$ as the infinum of a family of affine functions of $r$.

LEMMA 3.1. *One has*
$$\mathrm{RiskA}(\varepsilon) \leq \Phi_*(\ln(2/\varepsilon)).$$

PROOF. Given $\delta > 0$ and $\varepsilon \in (0, 1/4)$, let us build an affine estimate with $\varepsilon$-risk $\leq R \equiv \Phi_*(\ln(2/\varepsilon)) + \delta/2$, namely, as follows. By (3.1), there exist $\phi_* \in \mathcal{F}$ and $\alpha_* > 0$, such that

$$2\Phi_*(\ln(2/\varepsilon)) + \delta/2$$

$$\geq \max_{x,y \in X} \Phi_{\varepsilon/2}(x, y; \phi_*, \alpha_*)$$

$$= \underbrace{\max_{x \in X}\left[g^T x + \alpha_* \ln\left(\int_\Omega \exp\{-\alpha_*^{-1}\phi_*(\omega)\} p_{A(x)}(\omega) P(d\omega)\right) + \alpha_* \ln(2/\varepsilon)\right]}_{U}$$

$$\times \underbrace{\max_{y \in X}\left[-g^T y + \alpha_* \ln\left(\int_\Omega \exp\{\alpha_*^{-1}\phi_*(\omega)\} p_{A(y)}(\omega) P(d\omega)\right) + \alpha_* \ln(2/\varepsilon)\right]}_{V}.$$

Setting $c = \frac{U-V}{2}$, we have

$$\max_{x \in X}\left[g^T x + \alpha_* \ln\left(\int_\Omega \exp\{-\alpha_*^{-1}[\phi_*(\omega) + c]\} p_{A(x)}(\omega) P(d\omega)\right) + \alpha_* \ln(2/\varepsilon)\right]$$

$$= U - c = \frac{U+V}{2} \leq \Phi_*(\ln(2/\varepsilon)) + \delta/4 = R - \delta/4,$$

$$\max_{y \in Y}\left[g^T x + \alpha_* \ln\left(\int_\Omega \exp\{\alpha_*^{-1}[\phi_*(\omega) + c]\} p_{A(y)}(\omega) P(d\omega)\right) + \alpha_* \ln(2/\varepsilon)\right]$$

$$= V + c = \frac{U+V}{2} \leq \Phi_*(\ln(2/\varepsilon)) + \delta/4 = R - \delta/4$$

or, equivalently,

$$\max_{x \in X} \ln\left(\int_\Omega \exp\{\alpha_*^{-1}[g^T x - (\phi_*(\omega) + c) - R]\} p_{A(x)}(\omega) P(d\omega)\right)$$

$$\leq \ln(\varepsilon/2) - \frac{\delta}{4\alpha_*} \equiv \ln(\varepsilon'/2),$$

$$\max_{y \in X} \ln\left(\int_\Omega \exp\{\alpha_*^{-1}[(\phi_*(\omega) + c) - R - g^T y]\} p_{A(y)}(\omega) P(d\omega)\right) \leq \ln(\varepsilon'/2),$$

that is,

(a) $\forall x \in X$: $\quad \int_\Omega \exp\{\alpha_*^{-1}[g^T x - (\phi_*(\omega) + c) - R]\} p_{A(x)}(\omega) P(d\omega) \leq \varepsilon'/2,$

(b) $\forall y \in X$: $\quad \int_\Omega \exp\{\alpha_*^{-1}[[\phi_*(\omega) + c] - R - g^T y]\} p_{A(y)}(\omega) P(d\omega) \leq \varepsilon'/2.$

For a given $x \in X$, the exponent in (a) is nonnegative and is $> 1$, for all $\omega$ such that $g^T x - [\phi_*(\omega) + c] > R$; therefore, (a) implies that $\text{Prob}_{\omega \sim p_{A(x)}(\cdot)}\{g^T x > [\phi_*(\omega) + c] + R\} \leq \varepsilon'/2$, for every $x \in X$. By similar reasons, (b) implies that

$\text{Prob}_{\omega \sim p_{A(x)}(\cdot)}\{g^T x < [\phi_*(\omega) + c] - R\} \leq \varepsilon'/2$, for all $x \in X$. Since by construction $\varepsilon' < \varepsilon$, we see that the $\varepsilon$-risk of the affine estimate $\hat{g}(\omega) = \phi_*(\omega) + c$ is $\leq R$, as claimed. □

LEMMA 3.2. *One has*

(3.2) $\quad\quad\quad\quad \delta \in (0,1) \quad \Rightarrow \quad \text{Risk}_*(\delta^2/4) \geq \Phi_*(\ln(1/\delta)),$

*whence also*

(3.3) $\quad \varepsilon \in (0, 1/4) \quad \Rightarrow \quad \text{Risk}_*(\delta) \geq \dfrac{\ln(1/(4\varepsilon))}{2\ln(2/\varepsilon)} \Phi_*(\ln(2/\varepsilon)).$

PROOF. To prove (3.2), let us set $\rho = \ln(1/\delta)$. The function $\Psi_\rho(x,y) = \inf_{\phi \in \mathcal{F}, \alpha > 0} \Phi_\rho(x,y;\phi,\alpha)$ takes values in $\{-\infty\} \cup \mathbb{R}$, is upper semicontinuous (since $\Phi_r$ is continuous) and is not identically $-\infty$ (in fact, it is even $\geq 0$ when $y = x$). Thus, $\Psi_\rho$ achieves its maximum on $X \times X$ at certain point $(\bar{x}, \bar{y})$, and for any $(\alpha > 0, \phi \in \mathcal{F})$:

(3.4) $\quad \Phi_\rho(\bar{x}, \bar{y}; \phi, \alpha) \geq \Psi_\rho(\bar{x}, \bar{y}) = \sup_{x,y \in X} \inf_{\phi \in \mathcal{F}, \alpha > 0} \Phi_\rho(x,y;\phi,\alpha) = 2\Phi_*(\rho),$

where the concluding inequality is given by (3.1). Since $(\mathcal{D}, \mathcal{F})$ is a good pair, setting $\mu = A(\bar{x})$, $\nu = A(\bar{y})$ and $\bar{\phi}(\omega) = \frac{1}{2}\ln(p_\mu(\omega)/p_\nu(\omega))$, we get $\bar{\phi} \in \mathcal{F}$, which combines with (3.4) to imply that

$\forall (\alpha > 0):$

$$2\Phi_*(\rho) \leq g^T\bar{x} - g^T\bar{y} + \alpha\left[\ln\left(\int_\Omega \exp\{-\alpha^{-1}[\alpha\bar{\phi}(\omega)]\} p_\mu(\omega) P(d\omega)\right)\right.$$
$$\left. + \ln\left(\int_\Omega \exp\{\alpha^{-1}[\alpha\bar{\phi}(\omega)]\} p_\nu(\omega) P(d\omega)\right) + 2\rho\right]$$
$$= g^T\bar{x} - g^T\bar{y} + 2\alpha\left[\ln\left(\int_\Omega \sqrt{p_\mu(\omega)p_\nu(\omega)} P(d\omega)\right) + \rho\right].$$

The resulting inequality holds true for all $\alpha > 0$, meaning that

(3.5) $\quad\quad$ (a) $\quad g^T\bar{x} - g^T\bar{y} \geq 2\Phi_*(\rho) = 2\Phi_*(\ln(1/\delta)),$

$\quad\quad\quad\quad$ (b) $\quad \int_\Omega \sqrt{p_\mu(\omega) p_\nu(\omega)} P(d\omega) \geq \exp\{-\rho\} = \delta.$

Now assume, in contrast to what should be proved, that $\text{Risk}_*(\delta^2/4) < \Phi_*(\ln(1/\delta))$. Then, there exists $R' < \Phi_*(\ln(1/\delta))$, $\delta' < \delta^2/4$ and an estimate $\hat{g}(\omega)$ such that

$$\text{Prob}_{\omega \sim p_{A(x)}(\cdot)}\{|\hat{g}(\omega) - g^T x| > R'\} \leq \delta' \quad \forall x \in X.$$

Now, consider two hypotheses $\Pi_{1,2}$ on the distribution of $\omega$ stating that the densities of the distribution with regard to $P$ are $p_\mu$ and $p_\nu$, respectively. Consider a procedure for distinguishing between the hypotheses as follows: after $\omega$ is observed, we compare $\hat{g}(\omega)$ with $\bar{g} = \frac{1}{2}[g^T \bar{x} + g^T \bar{y}]$; if $\hat{g}(\omega) \geq \bar{g}$, we accept $\Pi_1$, otherwise we accept $\Pi_2$. Note that by (3.5)(a) and due to $R' < \Phi_*(\ln(1/\delta))$, the probability to accept $\Pi_2$ when $\Pi_1$ is true is $\leq$ the probability for $\hat{g}(\omega)$ to deviate from $g^T \bar{x}$ by at most $R'$, that is, it is $\leq \delta'$. Similarly, the probability to accept $\Pi_1$ when $\Pi_2$ is true is $\leq \delta'$. Now, let $\Omega_1$ be the part of $\Omega$ where our hypotheses testing routine accepts $\Pi_1$, so that in $\Omega_2 = \Omega \setminus \Omega_1$ the routine accepts $\Pi_2$. As we just have seen,

$$\int_{\Omega_1} p_\nu(\omega) P(d\omega) \leq \delta', \qquad \int_{\Omega_2} p_\mu(\omega) P(d\omega) \leq \delta',$$

whence

$$\int_\Omega \sqrt{p_\mu(\omega) p_\nu(\omega)} P(d\omega) = \sum_{i=1}^2 \int_{\Omega_i} \sqrt{p_\mu(\omega) p_\nu(\omega)} P(d\omega)$$
$$\leq \sum_{i=1}^2 \left(\int_{\Omega_i} p_\mu(\omega) P(d\omega)\right)^{1/2} \left(\int_{\Omega_i} p_\nu(\omega) P(d\omega)\right)^{1/2}$$
$$\leq 2\sqrt{\delta'} < 2\sqrt{\delta^2/4} = \delta.$$

The resulting inequality $\int_\Omega \sqrt{p_\mu(\omega) p_\nu(\omega)} P(d\omega) < \delta$ contradicts (3.5)(b); we have arrived at a desired contradiction. (3.2) is proved.

To prove (3.3), let us set $\delta = 2\sqrt{\varepsilon}$, so that $\text{Risk}_*(\varepsilon) = \text{Risk}_*(\delta^2/4) \geq \Phi_*(\ln(1/\delta)) = \Phi_*(\frac{1}{2}\ln(\frac{1}{4\varepsilon}))$, where the concluding $\geq$ is due to (3.2). Now recall that $\Phi_*(r)$ is a nonnegative and concave function of $r \geq 0$, so that $\Phi_*(tr) \geq t\Phi_*(r)$, for all $r \geq 0$ and $0 \leq t \leq 1$. We therefore have

$$\Phi_*\left(\frac{1}{2}\ln\left(\frac{1}{4\varepsilon}\right)\right) \geq \frac{\ln(1/(4\varepsilon))}{2\ln(2/\varepsilon)} \Phi_*\left(\ln\left(\frac{2}{\varepsilon}\right)\right)$$

and we arrive at (3.3). □

Lemmas 3.1 and 3.2 clearly imply Theorem 3.1. □

REMARK 3.1. Lemmas 3.1 and 3.2 provide certain information even beyond the case when the pair $(\mathcal{D}, \mathcal{F})$ is good, specifically, that:

(i) The $\varepsilon$-risk of an affine estimate can be made arbitrarily close to the quantity

$$\Phi_+(\varepsilon) = \inf_{\phi \in \mathcal{F}, \alpha > 0} \sup_{x,y \in X} \Phi_{\ln(2/\varepsilon)}(x, y; \phi, \alpha)$$

(cf. Lemma 3.1);

(ii) We have $\text{Risk}_*(\varepsilon) \geq \Phi_-(\varepsilon) = \sup_{x,y \in X} \inf_{\phi \in \mathcal{F}, \alpha > 0} \Phi_{1/2\ln(1/(4\varepsilon))}(x, y; \phi, \alpha)$ (cf. Lemma 3.2).

As it is seen from the proofs of Lemmas 3.1 and 3.2, both these statements hold true without the goodness assumption. The role of the latter is in ensuring that $\Phi_+(\varepsilon)$ is within an absolute constant factor of $\Phi_-(\varepsilon)$.

Lemma 3.2 Implies the following result.

PROPOSITION 3.1. *Under the premise of Theorem 3.1, the Hellinger affinity*

$$\text{AffH}(\mu, \nu) = \int_\Omega \sqrt{p_\mu(\omega) p_\nu(\omega)} P(d\omega)$$

*is a continuous and log-concave function on $\mathcal{M} \times \mathcal{M}$, and the quantity $\Phi_*(r)$, $r \geq 0$, admits the following representation:*

(3.6) $\quad 2\Phi_*(r) = \max_{x,y} \{g^T x - g^T y : \text{AffH}(A(x), A(y)) \geq \exp\{-r\}, x, y \in X\}.$

We see that the upper bound $\Phi_*(\ln(2/\varepsilon))$ on $\text{RiskAff}(\varepsilon)$ stated in Theorem 3.1 admits a very transparent interpretation: this bound is the maximum of the variation $\frac{1}{2} \max_{x,y} [g^T x - g^T y]$ of the estimated functional on the set of pairs $x, y \in X$ with the associated distributions "close" to each other, namely, such that $\text{AffH}(A(x), A(y)) \geq \varepsilon/2$. Observe that asymptotically (when $r$ becomes small),[2] $\Phi_*(r)$ is equivalent to the *modulus of continuity $\omega(r, X)$ of $g$ with regard to the Hellinger distance*, introduced in [7].

PROOF OF PROPOSITION 3.1. By exactly the same argument as in the proof of Theorem 3.1, the function $\Psi(\mu, \nu; \phi) : (\mathcal{M} \times \mathcal{M}) \times \mathcal{F} \to \mathbb{R}$,

$$\Psi(\mu, \nu; \phi) = \left[ \ln\left( \int \exp\{-\phi(\omega)\} p_\mu(\omega) P(d\omega) \right) \right.$$
$$\left. + \ln\left( \int \exp\{\phi(\omega)\} p_\nu(\omega) P(d\omega) \right) \right]$$

is well defined and continuous on its domain, and this function is convex in $\phi$ and concave in $(\mu, \nu)$. We claim that

(3.7) $\quad \ln(\text{AffH}(\mu, \nu)) = \frac{1}{2} \min_\phi \Psi(\mu, \nu; \psi),$

which would imply that $\ln(\text{AffH}(\cdot))$ is indeed a finite concave function on $\mathcal{M} \times \mathcal{M}$ and as such is continuous (recall that $\mathcal{M}$ is open). To justify our claim, note that, for fixed $\mu, \nu \in \mathcal{M}$, setting $\phi = \frac{1}{2}\ln(p_\nu/p_\mu)$, we get a function from $\mathcal{F}$ such that $\Psi(\mu, \nu; \bar{\phi}) = 2\ln(\text{AffH}(\mu, \nu))$. To complete the verification of (3.7), it suffices to demonstrate that $\Psi(\mu, \nu; \phi) \geq \Psi(\mu, \nu; \bar{\phi})$ whenever $\phi \in \mathcal{F}$, which is immediate, since setting $\phi = \bar{\phi} + \Delta$, we have

$$\exp\{\Psi(\mu, \nu; \bar{\phi})/2\}$$
$$= \int \sqrt{p_\mu(\omega) p_\nu(\omega)} P(d\omega)$$
$$= \int [(p_\mu(\omega) p_\nu(\omega))^{1/4} \exp\{-\Delta(\omega)/2\}]$$
$$\quad \times [(p_\mu(\omega) p_\nu(\omega))^{1/4} \exp\{\Delta(\omega)/2\}] P(d\omega)$$
$$\leq \left[\int \sqrt{p_\mu(\omega) p_\nu(\omega)} \exp\{-\Delta(\omega)\} P(d\omega)\right]^{1/2}$$
$$\quad \times \left[\int \sqrt{p_\mu(\omega) p_\nu(\omega)} \exp\{\Delta(\omega)\} P(d\omega)\right]^{1/2}$$
$$= \exp\{\Psi(\mu, \nu; \phi)/2\}.$$

---
[2]Recall that we consider here the case of one observation.

Now, note that by (3.1)

$$2\Phi_*(r) = \sup_{x,y \in X} \left\{ \inf_{\phi \in \mathcal{F}, \alpha > 0} [g^T x - g^T y + \alpha \Psi(A(x), A(y); \alpha^{-1}\phi) + 2\alpha r] \right\}$$
$$= \sup_{x,y \in X} \left\{ g^T x - g^T y + \inf_{\alpha > 0} \alpha \left[ \inf_{\phi \in \mathcal{F}} \Psi(A(x), A(y); \alpha^{-1}\phi) + 2r \right] \right\}$$
$$= \sup_{x,y \in X} \left\{ g^T x - g^T y + \inf_{\alpha > 0} \alpha \left[ \inf_{\psi \equiv \alpha^{-1}\phi \in \mathcal{F}} \Psi(A(x), A(y); \psi) + 2r \right] \right\}$$
$$= \sup_{x,y \in X} \left\{ g^T x - g^T y + \underbrace{\inf_{\alpha > 0} \alpha [2\ln(\text{AffH}(A(x), A(y))) + 2r]}_{= \begin{cases} 0, & \ln(\text{AffH}(A(x), A(y))) + r \geq 0, \\ -\infty, & \ln(\text{AffH}(A(x), A(y))) + r < 0 \end{cases}} \right\}$$

[see (3.7)]

$$= \max_{x,y} \{g^T x - g^T y : \text{AffH}(A(x), A(y)) \geq \exp\{-r\}, x, y \in X\}. \quad \square$$

3.2. *The case of multiple observations.* In Problem I, our goal was to estimate $g^T x$ from a *single* observation $\omega$ of the random variable $\omega \sim p_{A(x)}(\cdot)$, associated with $x$. The result can be immediately extended to the case when we want to recover $g^T x$ from a sample of independent observations $\omega_1, \ldots, \omega_L$ of random variables $\omega_\ell$ with distributions parameterized by $x$. Specifically, let $(\Omega_\ell, P_\ell)$ and $(\mathcal{D}_\ell, \mathcal{F}_\ell)$, $1 \leq \ell \leq L$, be as in Example 4, and let every pair $(\mathcal{D}_\ell, \mathcal{F}_\ell)$ be good. Assume, further, that $X \subset \mathbb{R}^n$ is a convex compact set and $A_\ell(x)$ are affine mappings with $A_\ell(X) \subset \mathcal{M}_\ell$. Given a linear form $g^T z$ on $\mathbb{R}^n$ and a sequence of independent realizations $\omega_\ell \sim p^\ell_{A_\ell(x)}(\cdot)$, $\ell = 1, \ldots, L$, we want to recover from these observations the value $g^T x$ of the given affine form at the "signal" $x$ underlying our observations.

In our current situation, we call a candidate estimate $\hat{g}(\omega_1, \ldots, \omega_L)$ affine if it is of the form

(3.8) $$\hat{g}(\omega_1, \ldots, \omega_L) = \sum_{\ell=1}^L \phi_\ell(\omega_\ell),$$

where $\phi_\ell \in \mathcal{F}_\ell$, $\ell = 1, \ldots, L$. Note that setting $(\mathcal{D}, \mathcal{F}) = \bigotimes_{\ell=1}^L (\mathcal{D}_\ell, \mathcal{F}_\ell)$, we reduce the situation to the one we have already considered. In particular, Theorem 3.1 along with the proof of Lemma 3.1 implies the following result (where the $\varepsilon$-risks—of an estimate, the minimax optimal and the affine-minimax optimal—are defined exactly as in the single-observation case).

THEOREM 3.2. *In the situation just described, for $r > 0$, let*

$$\Phi_r(x, y; \phi, \alpha) = \alpha \left[ \sum_{\ell=1}^L \ln \left( \int_{\Omega_\ell} \exp\{-\alpha^{-1}\phi_\ell(\omega_\ell)\} p^\ell_{A_\ell(x)}(\omega_\ell) P(d\omega_\ell) \right) \right.$$
$$\left. + \left( \int_{\Omega_\ell} \exp\{\alpha^{-1}\phi_\ell(\omega_\ell)\} p^\ell_{A_\ell(y)}(\omega_\ell) P(d\omega_\ell) \right) \right]$$
$$+ g^T x - g^T y + 2\alpha r: \qquad Z \times \mathcal{F}_+ \to \mathbb{R},$$
$$Z = X \times X,$$
$$\mathcal{F}_+ = \mathcal{F}_1 \times \cdots \times \mathcal{F}_L \times \{\alpha > 0\}.$$

*The function $\Phi_r$ is continuous on its domain, concave in the $(x, y)$-argument, convex in the $(\phi, \alpha)$-argument and possesses a well-defined saddle point value*

$$2\Phi_*(r) = \sup_{x,y \in X} \underbrace{\inf_{\phi, \alpha \in \mathcal{F}_+} \Phi_r(x, y; \phi, \alpha)}_{\underline{\Phi}_r(x,y)} = \inf_{(\phi, \alpha) \in \mathcal{F}_+} \underbrace{\sup_{x, y \in X} \Phi_r(x, y; \phi, \alpha)}_{\overline{\Phi}_r(\phi, \alpha)},$$

*which is a concave and nonnegative function of $r \geq 0$. Moreover:*

(i) *For all $\varepsilon \in (0, 1/4)$, we have*

$$\text{RiskA}(\varepsilon) \leq \Phi_*(\ln(2/\varepsilon)) \leq \theta(\varepsilon) \text{Risk}_*(\varepsilon), \qquad \theta(\varepsilon) = \frac{2\ln(2/\varepsilon)}{\ln(1/(4\varepsilon))}.$$

(ii) *Given $\varepsilon \in (0, 1/4)$ and $\delta \geq 0$, in order to build an affine estimate with $\varepsilon$-risk not exceeding $[\Phi_*(\ln(2/\varepsilon)) + \delta]$, where $\delta > 0$ is given, it suffices to find $\alpha_* > 0$ and $\phi_\ell^* \in \mathcal{F}_\ell$, $1 \leq \ell \leq L$, such that*

$$\overline{\Phi}_{\ln(2/\varepsilon)}(\phi^*, \alpha_*) \leq 2\Phi_*(\ln(2/\varepsilon)) + \delta/2,$$

*to compute the quantity*

$$c = \frac{1}{2} \max_{x \in X} \left[ g^T x + \alpha_* \sum_{\ell=1}^L \ln \left( \int_{\Omega_\ell} \exp\{-\alpha^{-1} \phi_\ell^*(\omega_\ell)\} p_{A_\ell(x)}^\ell(\omega_\ell) P_\ell(d\omega_\ell) \right) \right]$$
$$- \frac{1}{2} \max_{y \in X} \left[ -g^T y + \alpha_* \sum_{\ell=1}^L \ln \left( \int_{\Omega_\ell} \exp\{\alpha^{-1} \phi_\ell^*(\omega_\ell)\} p_{A_\ell(y)}^\ell(\omega_\ell) P_\ell(d\omega_\ell) \right) \right]$$

*and to set*

$$(3.9) \qquad \hat{g}(\omega_1, \ldots, \omega_L) = \sum_{\ell=1}^L \phi_\ell^*(\omega_\ell) + c.$$

REMARK 3.2. Computing the "nearly optimal" affine estimate (3.9) reduces to convex programming and thus can be carried out efficiently, provided that we are given explicit descriptions of:

- the linear spaces $\mathcal{F}_\ell$, $\ell = 1, \ldots, L$ (as it is the case, e.g., in Examples 1–3),
- and $X$ (e.g., by a list of efficiently computable convex constraints which cut $X$ out of $\mathbb{R}^n$) and are capable to compute efficiently the value of $\Phi_r$ at a point.

REMARK 3.3. Assume that the observations $\omega_\ell$, $\ell_0 \leq \ell \leq \ell_1$, are copies of the same random variable [i.e., $\Omega_\ell, P_\ell, \mathcal{D}_\ell, \mathcal{F}_\ell, A_\ell(\cdot)$ are independent of $\ell$ for $\ell_0 \leq \ell \leq \ell_1$]. Then, the convex function $\overline{\Phi}_r(\phi_1, \ldots, \phi_L, \alpha)$ is symmetric with regard to the arguments $\phi_\ell \in \mathcal{F}_{\ell_0}$, $\ell_0 \leq \ell \leq \ell_1$, and therefore, when building the estimate (3.9) we lose nothing when restricting ourselves to $\phi$'s satisfying $\phi_\ell = \phi_{\ell_0}$, $\ell_0 \leq \ell \leq \ell_1$, which allows to reduce the computational effort of building $\alpha_*, \phi_\ell^*$.

3.2.1. *Illustration.* Consider the toy problem where we want to recover the probability $p$ of getting 1 from a Bernoulli distribution, given $L$ independent realizations $\omega_1, \ldots, \omega_L$ of the associated random variable. To handle the problem, we specialize our general setup as follows:

- $(\Omega_\ell, P_\ell)$, $1 \leq \ell \leq L$, are identical to the two-point set $\{0; 1\}$ with the counting measure;
- $\mathcal{M}$ is the interval $(0, 1)$, and $p_\mu(1) = 1 - p_\mu(0) = \mu$, $\mu \in \mathcal{M}$;
- $X$ is a compact convex subset in $\mathcal{M}$, say, the segment $[1\cdot\text{e--}16, 1\text{--}1/\text{e--}16]$, and $A(x) = x$.

TABLE 1
*Recovering the parameter of a Bernoulli distribution*

| $\varepsilon$ | $L$ | $\gamma$ | $\delta$ | Upper risk bound | Lower risk bound | Ratio of bounds | $\vartheta(\varepsilon)$ |
|---|---|---|---|---|---|---|---|
| 0.05 | 10 | 2.91e–1 | 4.18e–2 | 3.61e–1 | 2.49e–1 | 1.45 | 4.58 |
| 0.05 | 100 | 4.13e–2 | 9.17e–3 | 1.33e–1 | 8.19e–2 | 1.63 | 4.58 |
| 0.05 | 1000 | 4.29e–3 | 9.91e–4 | 4.29e–3 | 2.60e–3 | 1.65 | 4.58 |
| 0.01 | 10 | 3.58e–1 | 2.83e–2 | 4.04e–1 | 3.29e–1 | 1.23 | 3.29 |
| 0.01 | 100 | 5.83e–2 | 8.84e–2 | 1.59e–1 | 1.15e–1 | 1.38 | 3.29 |
| 0.01 | 1000 | 6.15e–3 | 9.88e–4 | 5.13e–2 | 3.67e–3 | 1.40 | 3.29 |
| 0.001 | 10 | 4.19e–1 | 1.61e–2 | 4.42e–1 | 3.98e–1 | 1.11 | 2.75 |
| 0.001 | 100 | 8.15e–2 | 8.37e–3 | 1.88e–1 | 1.51e–1 | 1.24 | 2.75 |
| 0.001 | 1000 | 8.79e–3 | 9.82e–4 | 6.14e–3 | 4.88e–3 | 1.26 | 2.75 |

Invoking Remark 3.3, we lose nothing when restricting ourselves to affine estimates of the form (3.8) with mutually identical functions $\phi_\ell(\cdot)$, $1 \leq \ell \leq L$, that is, with the estimates

$$\hat{g}(\omega_1, \ldots, \omega_L) = \gamma + \delta \sum_{\ell=1}^{L} \omega_\ell.$$

Invoking Theorem 3.2, the coefficients $\gamma$ and $\delta$ are readily given by the $\phi$-component of the saddle point (max in $x, y \in X$, min in $\phi = [\phi_0; \phi_1] \in \mathbb{R}^2$ and $\alpha > 0$) of the convex–concave function

$$x - y + \alpha[L \ln(\varepsilon^{-\phi_0/\alpha}(1-x) + \varepsilon^{-\phi_1/\alpha}x) \\ + L \ln(\varepsilon^{\phi_0/\alpha}(1-y) + \varepsilon^{\phi_1/\alpha}y) + 2\ln(2/\varepsilon)];$$

the (guaranteed upper bound on the) $\varepsilon$-risk of this estimate is half of the corresponding saddle point value. The saddle point (it is easily seen that it does exist) can be computed with high accuracy by standard convex programming techniques. In Table 1, we present the nearly optimal affine estimates along with the corresponding risks. In the table, the upper risk bound is the one guaranteed by Theorem 3.2 and the lower risk bound is the largest $d$ such that the hypotheses "$p = 0.5 + d$" and "$p = 0.5 - d$" cannot be distinguished from $L$ independent observations of a random variable $\sim$ Bernoulli($p$) with the sum of probabilities of errors $< 2\varepsilon$ [this easily computable quantity is a lower bound on the minimax optimal $\varepsilon$-risk Risk$_*(\varepsilon)$], and $\vartheta(\varepsilon) = \frac{2\ln(2/\varepsilon)}{\ln(0.25/\varepsilon)}$ is the theoretical upper bound on the "level of nonoptimality" of our estimate. As it could be guessed in advance, for large $L$, the near-optimal affine estimate is close to the trivial estimate $\frac{1}{L}\sum_{\ell=1}^{L}\omega_\ell$.

**4. Applications.** In this section, we present some applications of Theorems 3.1 and 3.2.

4.1. *Positron emission tomography.* The positron emission tomography (PET) is a noninvasive diagnostic tool allowing us to visualize not only the anatomy of tissues in a body, but their functioning as well. In PET, a patient is administered a radioactive tracer chosen in such a way that it concentrates in the areas of interest (e.g., those of high metabolic activity in early diagnosis of cancer). The tracer disintegrates, emitting positrons which then annihilate with nearby electrons to produce pairs of photons flying at the speed of light in opposite directions; the orientation of the resulting *line of response* (LOR) is completely random. The patient is placed in a cylinder with the surface split into small detector cells. When two of the detectors are hit by photons "nearly simultaneously"—within an appropriately chosen short time window—the event indicates that somewhere at a line crossing the detectors a disintegration act took place. Such an event is registered, and the data collected by the PET device form a list of the number of events registered in every one of the *bins* (pairs of detectors) in the course of a given time $t$. The goal of a PET reconstruction algorithm is to recover the density of the tracer from this data. The standard mathematical model of PET is as follows. After discretization of the field of view, there are $N$ voxels (small 3D cubes) assigned with nonnegative (and unknown) amounts $x_i$ of the traces $i = 1, \ldots, n$. The number of LORs emanating from a voxel $i$ is a realization of a Poisson random variable with parameter $x_i$, and these variables for different voxels are independent. Every LOR emanating from a voxel $i$ is subject to a "lottery," which decides in which bin (pair of detectors) it will be registered or if it will be registered at all—some LORs can intersect the surface of the cylinder only in one point or not intersect it at all and thus are missed. The role of the lottery is played by the random orientation of the LOR in question, and outcomes of different lotteries are independent. The probabilities $q_{i\ell}$ for a LOR emanating from voxel $i$ to be registered in bin $\ell$ are known (they are readily given by the geometry of the device). With this model, the data registered by PET is a realization of a random vector $(\omega_1, \ldots, \omega_L)$ ($L$ is the total number of bins) with independent Poisson-distributed coordinates, the parameter of the Poisson distribution associated with $\omega_\ell$ being

$$A_\ell(x) = \sum_{i=1}^n q_{i\ell} x_i.$$

Assume that our a priori information on $x$ allows us to point out a convex compact set $X \subset \{x \in \mathbb{R}^n : x > 0\}$, such that $x \in X$. Assuming without loss of generality that $\sum_i q_{i\ell} > 0$ for every $\ell$ (indeed, we can eliminate all bins $\ell$ which never register LORs) and invoking Example 2, we find ourselves in the situation of Section 3.2. It follows that in order to evaluate a given linear form $g^T x$ of the unknown tracer density $x$, we can use the construction from Theorem 3.2 to build a near-optimal affine estimate of $g^T x$. The recipe suggested to this end by Theorem 3.2 reads as follows: the estimate is of the form

$$\hat{g}(\omega) = \sum_{\ell=1}^{L} \gamma_\ell^* y_\ell + c_*,$$

where $y_\ell$ is the number of LORs registered in bin $\ell$ and $\gamma^* = [\gamma_1^*; \ldots; \gamma_L^*]$, $c_*$ are given by an optimal solution $(\gamma^*, \alpha_*)$ to the convex optimization problem

$$\min_{\alpha > 0, \gamma} \overline{\Phi}_r(\gamma, \alpha),$$

$$\overline{\Phi}_r(\gamma, \alpha) = \max_{x, y \in X} \Bigg\{ g^T x - g^T y $$
$$+ \alpha \Bigg[ \sum_{\ell=1}^{L} [q_\ell(x) \exp\{-\alpha^{-1}\gamma_\ell\} + q_\ell(y) \exp\{\alpha^{-1}\gamma_\ell\}] $$
$$- q(x) - q(y) + 2r \Bigg] \Bigg\},$$

$$r = \ln(2/\varepsilon), \qquad q_\ell(z) = \sum_{i=1}^{n} q_{i\ell} z_i, \qquad q(z) = \sum_{\ell=1}^{L} q_\ell(z).$$

It is easily seen that the problem is solvable with

$$c_* = \frac{1}{2} \Bigg[ \max_{x \in X} \Bigg\{ g^T x + \alpha_* \Bigg[ -q(x) + \sum_{\ell=1}^{L} q_\ell(x) \exp\{-\alpha_*^{-1}\gamma_\ell^*\} \Bigg] \Bigg\} $$
$$- \max_{y \in X} \Bigg\{ -g^T y + \alpha_* \Bigg[ -q(y) + \sum_{\ell=1}^{L} q_\ell(y) \exp\{\alpha_*^{-1}\gamma_\ell^*\} \Bigg] \Bigg\} \Bigg].$$

4.2. *Gaussian observations.* Now consider the standard problem of recovering a linear form $g^T x$ of a signal $x$ known to belong to a given convex compact set $X \subset \mathbb{R}^n$ via indirect observations of the signal corrupted by Gaussian noise. Without loss of generality, let the model of observations be

(4.1) $$\omega = Ax + \xi, \qquad \xi \sim \mathcal{N}(0, I_L).$$

The associated pair $(\mathcal{D}, \mathcal{F})$ is comprised of the shifts of the standard Gaussian distribution ($\mathcal{D}$) and all affine forms on $\mathbb{R}^L$ ($\mathcal{F}$) and is good (see Example 3). The affine estimates in the case in question are just the affine functions of $\omega$. The near-optimality of affine estimates in the case in question was established by Donoho [5], not only for the $\varepsilon$-risk, but for all risks based on the standard loss functions. We have the following direct corollary of Theorem 3.2 (cf. Theorem 2 and Corollary 1 of [5]):

PROPOSITION 4.1. *In the situation in question, the affine estimate $\hat{g}_\varepsilon(\cdot)$ yielded by Theorem 3.2 is asymptotically ($\varepsilon \to +0$) optimal, specifically,*

$$\varepsilon \in (0, 1/2) \quad \Rightarrow \quad \mathrm{Risk}(\hat{g}_\varepsilon; \varepsilon) \leq \psi(\varepsilon)\, \mathrm{Risk}_*(\varepsilon),$$

$$\psi(\varepsilon) = \frac{\sqrt{2\ln(2/\varepsilon)}}{\mathrm{ErfInv}(\varepsilon)} = 1 + o(1) \quad \text{as } \varepsilon \to +0$$

*[here $x = \mathrm{ErfInv}(y)$ stands for the inverse error function, i.e., $y = \frac{1}{\sqrt{2\pi}} \times \int_x^\infty e^{-t^2/2}\, dt$].*

PROOF. Let $G(\cdot)$ be the density of the $\mathcal{N}(0, I_L)$ distribution. By Theorem 3.2, we have $\mathrm{Risk}(\hat{g}_\varepsilon; \varepsilon) \leq \Phi_*(\ln(2/\varepsilon))$, where, for $r > 0$,

$$2\Phi_*(r) = \max_{x, y \in X} \underline{\Phi}_r(x, y),$$

$$\underline{\Phi}_r(x, y) = \inf_{\phi \in \mathbb{R}^L, \alpha > 0} \Big\{ g^T x - g^T y$$

$$+ \alpha \Big[ \ln\Big( \int \exp\{-\alpha^{-1} \phi^T \omega\} G(\omega - Ax)\, d\omega \Big)$$

$$+ \ln\Big( \int \exp\{\alpha^{-1} \phi^T \omega\} G(\omega - Ay)\, d\omega \Big) + 2r \Big] \Big\}$$

$$= \inf_{\phi \in \mathbb{R}^L, \alpha > 0} \Big\{ g^T x - g^T y + \phi^T A(y - x) + 2 \Big[ \alpha^{-1} \frac{\phi^T \phi}{2} + \alpha r \Big] \Big\}$$

$$= \inf_\phi \{ g^T x - g^T y + \phi^T A(x - y) + 2\sqrt{2r} \|\phi\|_2 \}$$

$$= \begin{cases} g^T x - g^T y, & \|A(x-y)\|_2 \leq 2\sqrt{2r}, \\ -\infty, & \|A(x-y)\|_2 > 2\sqrt{2r}. \end{cases}$$

Thus,

(4.2) $\quad \mathrm{Risk}(\hat{g}_\varepsilon; \varepsilon) \leq \Phi_*(\ln(2/\varepsilon)) = \frac{1}{2}[g^T \bar{x} - g^T \bar{y}]$

for certain $\bar{x}, \bar{y} \in X$ with $\|A(\bar{x}-\bar{y})\|_2 \leq 2\sqrt{2\ln(2/\varepsilon)}$. It remains to prove that

(4.3) $\quad \mathrm{Risk}_*(\varepsilon) \geq \psi^{-1}(\varepsilon)\tfrac{1}{2}\Phi_*(\ln(2/\varepsilon)).$

To this end, assume, on the contrary to what should be proved, that

$$\mathrm{Risk}_*(\varepsilon) < \psi^{-1}(\varepsilon) \Phi_*(\ln(2/\varepsilon)) \quad (= \psi^{-1}(\varepsilon) \tfrac{1}{2}[g^T \bar{x} - g^T \bar{y}]),$$

and let us lead this assumption to a contradiction. Under our assumption, there exists $\rho < \psi^{-1}(\varepsilon) \tfrac{1}{2}[g^T \bar{x} - g^T \bar{y}]$, $\varepsilon' < \varepsilon$ and an estimate $\widetilde{g}$ such that

(4.4) $\quad \forall (x \in X): \quad \mathrm{Prob}\{|\widetilde{g}(Ax + \xi) - g^T x| \geq \rho\} \leq \varepsilon'.$

Observing that $\psi(\varepsilon) > 1$, we see that $2\rho < [g^T \bar{x} - g^T \bar{y}]$. Let $\hat{x} = \bar{x}$ and $\hat{y}$ be a convex combination of $\bar{x}$ and $\bar{y}$ such that $2\rho = [g^T \hat{x} - g^T \hat{y}]$. Note that

$$\|A(\hat{x} - \hat{y})\|_2 = \underbrace{\left[\frac{2\rho}{[g^T \bar{x} - g^T \bar{y}]}\right]}_{<\psi^{-1}(\varepsilon)} \|A(\bar{x} - \bar{y})\|_2 \leq \psi^{-1}(\varepsilon) 2\sqrt{2\ln(2/\varepsilon)} = 2\,\text{erfinv}(\varepsilon).$$

Now, let $\Pi_1$ be the hypothesis that the distribution of an observation (4.1) comes from $x = \hat{x}$, and let $\Pi_2$ be the hypothesis that this distribution comes from $x = \hat{y}$. From (4.4) by the same standard argument as in the proof of Lemma 3.2, it follows that there exists a routine, based on a single observation (4.1), for distinguishing between $\Pi_1$ and $\Pi_2$, which rejects $\Pi_i$ when this hypothesis is true with probability $\leq \varepsilon'$, $i = 1, 2$. But, it is well known that the hypotheses on shifts of the standard Gaussian distribution indeed can be distinguished with the outlined reliability. This is possible if and only if the Euclidean distance between the corresponding shifts is at least $2\,\text{erfinv}(\varepsilon')$. This condition is *not* satisfied for our $\Pi_i$, $i = 1, 2$, which correspond to shifts $A\hat{x}$ and $A\hat{y}$, since $\|A\hat{x} - A\hat{y}\|_2 \leq 2\,\text{erfinv}(\varepsilon) < 2\,\text{erfinv}(\varepsilon)$. We have arrived at a desired contradiction. $\square$

In fact, the reasoning can be slightly simplified and strengthened to yield the following result.

PROPOSITION 4.2. *In the situation of Proposition 4.1, one can build efficiently an affine estimate $\hat{g}_\varepsilon$, such that*

$$0 < \varepsilon < 1/2 \quad \Rightarrow \quad \text{Risk}(\hat{g}_\varepsilon; \varepsilon) \leq \frac{\text{ErfInv}(\varepsilon/2)}{\text{ErfInv}(\varepsilon)} \text{Risk}_*(\varepsilon)$$

*[cf. Proposition 4.1, and note that $\frac{\text{ErfInv}(\varepsilon/2)}{\text{ErfInv}(\varepsilon)} < \frac{\sqrt{2\ln(2/\varepsilon)}}{\text{ErfInv}(\varepsilon)}$].*

PROOF. Let

$$\Psi(x, y; \phi) = g^T x - g^T y + \phi^T A(y - x) + 2\,\text{erfinv}(\varepsilon/2)\|\phi\|_2 : (X \times X) \times \mathbb{R}^L \to \mathbb{R}.$$

$\Psi$ clearly is a function which is continuous, convex in $\phi$ and concave in $(x, y)$ on its domain; by the same argument as in the proof of Theorem 3.1, $\Psi$ has a well-defined saddle point value

$$2\Psi_*(\varepsilon) = \inf_\phi \overbrace{\max_{x,y \in X} \Psi(x, y; \phi)}^{\overline{\Psi}(\phi)} = \max_{x,y \in X} \overbrace{\inf_\phi \Psi(x, y; \phi)}^{\underline{\Psi}(x,y)}.$$

The function

$$\overline{\Psi}(\phi) = \max_{x,y \in X}[g^T x - g^T y + \phi^T(Ay - Ax)] + 2\,\text{erfinv}(\varepsilon/2)\|\phi\|_2 \geq 2\,\text{erfinv}(\varepsilon)\|\phi\|_2$$

is a finite convex function on $\mathbb{R}^L$, which goes to $\infty$ as $\|\phi\|_2 \to \infty$, and therefore it attains its minimum at a point $\phi_*$, so that

$$2\Psi_*(\varepsilon) = \overline{\Psi}(\phi_*).$$

Setting

$$c_* = \frac{1}{2}\left[\max_{x \in X}[g^T x - \phi_*^T A x] - \max_{y \in Y}[-g^T y + \phi_*^T A y]\right],$$

we have, similar to the proof of Lemma 3.1, the following:

(a) $\quad \max_{x \in X}[g^T x - \phi_*^T A x - c_*] + \mathrm{erfinv}(\varepsilon/2)\|\phi_*\|_2 = \Psi_*(\varepsilon),$

(b) $\quad \max_{y \in X}[-g^T y + \phi_*^T A x + c_*] + \mathrm{erfinv}(\varepsilon/2)\|\phi_*\|_2 = \Psi_*(\varepsilon).$

Now, consider the affine estimate

$$\hat{g}_\varepsilon(\omega) = \phi_*^T \omega + c_*.$$

From (a) it follows that

$$\forall d > \Psi_*(\varepsilon): \quad \sup_{x \in X} \mathrm{Prob}\{g^T x - \hat{g}_\varepsilon(Ax + \xi) > d\} \le \varepsilon' < \varepsilon/2,$$

while (b) implies that

$$\forall d > \Psi_*(\varepsilon): \quad \sup_{y \in X} \mathrm{Prob}\{\hat{g}_\varepsilon(Ay + \xi) - g^T y > d\} \le \varepsilon' < \varepsilon/2.$$

We conclude that $\mathrm{Risk}(\hat{g}_\varepsilon; \varepsilon) \le \Psi_*(\varepsilon)$. To complete the proof, it suffices to demonstrate that

(4.5) $$\mathrm{Risk}_*(\varepsilon) \le \frac{\mathrm{ErfInv}(\varepsilon/2)}{\mathrm{ErfInv}(\varepsilon)}\Psi_*(\varepsilon).$$

To this end, observe that

$$\underline{\Psi}(x, y) = [g^T x - g^T y] + \inf_\phi\{\phi^T A(y - x) + 2\,\mathrm{erfinv}(\varepsilon/2)\|\phi\|_2\}$$

$$= \begin{cases} g^T x - g^T y, & \|A(y - x)\|_2 \le 2\,\mathrm{erfinv}(\varepsilon/2), \\ -\infty, & \text{otherwise,} \end{cases}$$

whence

$$\mathrm{Risk}_*(\hat{g}_\varepsilon; \varepsilon) \le \Psi_*(\varepsilon) = \tfrac{1}{2}[g^T \bar{x} - g^T \bar{y}],$$

for certain $\bar{x}, \bar{y} \in X$, such that $\|A(\bar{x} - \bar{y})\|_2 \le 2\,\mathrm{erfinv}(\varepsilon)$. Relation (4.5) can be derived from this observation by exactly the same argument as used in the proof of Proposition 4.1 to derive (4.3) from (4.2). □

**5. Adaptive version of the estimate.** In the situation of Problem I, let $X^1 \subset X^2 \subset \cdots \subset X^K$ be a nested collection of nonempty convex compact sets in $\mathbb{R}^n$, such that $A(X^K) \subset \mathcal{M}$. Consider a modification of the problem where the signal $x$ underlying our observation is known to belong to one of $X^k$ with value of $k \leq K$ unknown in advance. Given a linear form $g^T z$ on $\mathbb{R}^n$, let $\text{Risk}^k(\hat{g}; \varepsilon)$ and $\text{Risk}_*^k(\varepsilon)$ be, respectively, the $\varepsilon$-risk of an estimate $\hat{g}$ on $X^k$ and the minimax optimal $\varepsilon$-risk of recovering $g^T x$ on $X^k$. Let also $\Phi_*^k(r)$ be the function associated with $X = X^k$ according to (3.1). As it is immediately seen, the functions $\Phi_*^k(r)$ grow with $k$. Our goal is to modify the estimate $\hat{g}$ yielded by Theorem 3.1 in such a way that the $\varepsilon$-risk of the modified estimate on $X^k$ will be "nearly" $\text{Risk}_*^k(\varepsilon)$ *for every $k \leq K$*. This goal can be achieved by a straightforward application of the well-known Lepskii's adaptation scheme [19, 20] as follows.

Given $\delta > 0$, let $\delta' \in (0, \delta)$, and let $\hat{g}^k(\cdot)$ be the affine estimate with the $(\varepsilon/K)$-risk on $X^k$ not exceeding $\Phi_*^k(\ln(2K/\varepsilon)) + \delta'$ provided by Theorem 3.1 as applied with $\varepsilon/K$ substituted for $\varepsilon$ and $X^k$ substituted for $X$. Then, for any $k \leq K$,

$$\sup_{x \in X^k} \text{Prob}_{\omega \sim p_{A(x)}(\cdot)} \{|\hat{g}^k(\omega) - g^T x| > \Phi_*^k(\ln(2K/\varepsilon)) + \delta\}$$
(5.1)
$$\leq \varepsilon'/K < \varepsilon/K.$$

Given observation $\omega$, let us say that an index $k \leq K$ is $\omega$-*good*, if for any $k', k \leq k' \leq K$,

$$|\hat{g}^{k'}(\omega) - \hat{g}^k(\omega)| \leq \Phi_*^k(\ln(2K/\varepsilon)) + \Phi_*^{k'}(\ln(2K/\varepsilon)) + 2\delta.$$

Note that $\omega$-good indexes do exist (e.g., $k = K$). Given $\omega$, we can find the smallest $\omega$-good index $k = k(\omega)$; our estimate is nothing but $\hat{g}(\omega) = \hat{g}^{k(\omega)}(\omega)$.

PROPOSITION 5.1. *Assume that $\varepsilon \in (0, 1/4)$, and let*

$$\vartheta = 3 \frac{\ln(2K/\varepsilon)}{\ln(2/\varepsilon)}.$$

*Then, for any $(k, 1 \leq k \leq K)$,*

(5.2) $\quad \sup_{x \in X^k} \text{Prob}_{\omega \sim p_{A(x)}(\cdot)} \{|\hat{g}(\omega) - g^T x| > \vartheta \Phi_*^k(\ln(2K/\varepsilon)) + 3\delta\} < \varepsilon,$

*whence also*

(5.3) $\quad \forall (k, 1 \leq k \leq K): \quad \text{Risk}^k(\hat{g}; \varepsilon) \leq \frac{6 \ln((2K)/\varepsilon)}{\ln(1/(4\varepsilon))} \text{Risk}_*^k(\varepsilon) + 3\delta.$

PROOF. Setting $r = \ln(2K/\varepsilon)$, let us fix $\bar{k} \leq K$ and $x \in X^{\bar{k}}$ and call a realization $\omega$ *x-good*, if

(5.4) $\qquad \forall(k, \bar{k} \leq k \leq K): \qquad |\hat{g}^k(\omega) - g^T x| \leq \Phi_*^k(r) + \delta.$

Since $X^k \supset X^{\bar{k}}$ when $k \geq \bar{k}$, (5.1) implies that

$$\text{Prob}_{\omega \sim p_{A(x)}(\cdot)}\{\omega \text{ is good}\} \geq 1 - \varepsilon'.$$

Now, when $x$ is the signal and $\omega$ is $x$-good, relations (5.4) imply that $\bar{k}$ is an $\omega$-good index, so that $k(\omega) \leq \bar{k}$. Since $k(\omega)$ is an $\omega$-good index, we have

$$|\hat{g}(\omega) - \hat{g}^{\bar{k}}(\omega)| = |\hat{g}^{k(\omega)}(\omega) - \hat{g}^{\bar{k}}(\omega)| \leq \Phi_*^k(r) + \Phi_*^{\bar{k}}(r) + 2\delta,$$

which combines with (5.4) to imply that

(5.5) $\qquad |\hat{g}(\omega) - g^T x| \leq 2\Phi_*^{\bar{k}}(r) + \Phi_*^{k(\omega)}(r) + 3\delta \leq 3\Phi_*^{\bar{k}}(r) + 3\delta,$

where the concluding inequality is due to $k(\omega) \leq \bar{k}$ and to the fact that $\Phi_*^k$ grows with $k$. The bound (5.5) holds true whenever $\omega$ is $x$-good, which, as we have seen, happens with probability $\geq 1 - \varepsilon'$. Since $\varepsilon' < \varepsilon$ and $\bar{x} \in X^{\bar{k}}$ is arbitrary, we conclude that

(5.6) $\qquad \text{Risk}^{\bar{k}}(\hat{g}; \varepsilon) \leq 3\Phi_*^{\bar{k}}(r) + 3\delta.$

Using the nonnegativity and concavity of $\Phi_*^{\bar{k}}(\cdot)$ on the nonnegative ray and recalling the definition of $r$, we obtain $\Phi_*^k(r) \leq \frac{\ln(2K/\varepsilon)}{\ln(2/\varepsilon)}\Phi_*^k(\ln(2/\varepsilon))$ whenever $\varepsilon \leq 1/2$ and $K \geq 1$. Recalling the definition of $\vartheta$, the right-hand side in (5.6) therefore does not exceed $\vartheta \Phi_*^{\bar{k}}(\ln(2/\varepsilon)) + 3\delta$. Since $\bar{k} \leq K$ is arbitrary, we have proved (5.2). This bound, due to Lemma 3.2, implies (5.3). $\square$

**Acknowledgments.** The authors would like to acknowledge the valuable suggestions made by L. Birgè, Universitè Paris 6, and Alexander Goldenshluger, Haifa University.

Laboratoire Jean Kuntzmann  
Université Grenoble I  
51 rue des Mathématiques  
BP 53  
38041 Grenoble Cedex 9  
France  
E-mail: juditsky@imag.fr  

School of Industrial  
    and Systems Engineering  
Georgia Institute of Technology  
765 Ferst Drive  
Atlanta, Georgia 30332-0205  
USA  
E-mail: nemirovs@isye.gatech.edu